\newif\ifNoK
\newif\ifShort
\newif\ifOdsStyle

\NoKfalse 	

\Shorttrue

\OdsStyletrue 

\ifOdsStyle
\documentclass{svproc}
\else
\documentclass{article}
\fi

\usepackage{url}
\usepackage[normalem]{ulem}
\usepackage[leqno]{amsmath}
\usepackage{amsthm}
\usepackage{amssymb}
\usepackage{todonotes}
\usepackage{mathtools}
\usepackage{amsmath}
\usepackage{enumerate}
\usepackage{enumitem}
\usepackage{pdflscape}
\usepackage{algorithmic}

\newcommand{\unaryminus}{\scalebox{0.5}[1.0]{\( - \)}}
\newcommand{\enumWithLetter}[2]{$#1#2$}

\begin{document}

\ifOdsStyle

\mainmatter
\title{A Mixed-integer Linear Program to create the shifts in a supermarket}
\titlerunning{MILP for supermarket shifts}

\author{Nicol{\`o} Gusmeroli \and Andrea Bettinelli}
\authorrunning{Nicol{\`o} Gusmeroli et al.} 

\tocauthor{Nicol{\`o} Gusmeroli, and Andrea Bettinelli}
\institute{OPTIT srl, Bologna, Italy, \\
\email{nicolo.gusmeroli@optit.net, andrea.bettinelli@optit.net}}
\institute{OPTIT srl, Via Mazzini 82, 40135 Bologna, Italy,\\
\email{nicolo.gusmeroli@optit.net}, \email{andrea.bettinelli@optit.net} }

\else

\title{A Mixed-integer Linear Program to create the shifts in a supermarket}
\author{Nicol{\`o} Gusmeroli\thanks{OPTIT srl, Bologna, Italy. 
Emails: \{nicolo.gusmeroli, andrea.bettinelli@optit.net\} }
\and Andrea Bettinelli\footnotemark[1]}
\date{\today}
\fi

\maketitle

\begin{abstract}
The shift design and the personnel scheduling problem 
is known to be a difficult problem.
It is a real-world problem which has lots of applications
in the organization of 
\ifShort
companies.
\fi
Solutions are usually found by dividing the problem in two steps:
first the shifts are created,
then the employees are assigned to them 
by respecting a bunch of constraints.
The assignment of different tasks increases the complexity, 
since we have to consider the skills of the single employee 
necessary to perform any activity. 
In this paper we present a mixed-integer linear programming formulation
which models together the shift creation and 
the construction of rosters for employees,
with the objective of minimizing the amount of uncovered demand.
Finally we provide the results for
three real-world instances, 
confirming that this approach is promising.
\ifOdsStyle
\else
\fi
\end{abstract}

\section{Introduction}

Building employee rosters while respecting legal 
and organizational constraints to satisfy personnel requirements 
is an NP-complete problem, see~\cite{SMPTSP}.
In the past,  companies mainly adopted two techniques 
to simplify the problem:
\ifShort
one was to fix the shift hours beforehand
\fi
and then the employees were assigned in these shifts,
in the other they repeated the schedules in a cyclic way over the weeks
by shifting the employees, so the planning was done only once.
This second approach left the possibility that some worker was not available,
which was solved manually at any occurrence.
Nowadays the situation has changed since the employees 
want to have more freedom in the assignment of the shifts, 
hence the companies started to let people decide the rest day in the week 
and to give preferences on the working hours.
This flexibility made the mathematical modelling of this problem 
even more difficult and, most likely, less studied.

The problem studied consists in designing individual schedules 
for a group of heterogeneous workers, 
which amounts to fixing days-off, 
creating shifts, and assigning available fixed tasks 
within these shifts, by respecting both the governor and the company rules,
and it is called the 
Shift-Design Personnel Task Scheduling Problem (SDPTSP). 
Due to the difficulty of the problem the related literature work is scarce:
to the best of our knowledge only few approaches have been proposed 
and no paper ever presented 
an explicit mathematical formulation. 
We aim to mind this gap by formulating a 
mixed-integer linear programming formulation (MILP) 
for the SDPTSP.
The motivation for this work comes from a real-world application, 
namely the creation and the assignment of the shifts 
to the employees of a large Italian retail company,
hence we consider additional constraints for our specific use case.
The aim of this collaboration is threefold:
automatize the rostering process, reduce the work of the managers, and 
improve the current schedules.
Together with our customer we developed a model 
which was solved with standard algorithms, 
i.e., not specific for this problem, 
so there is big room for improvements.
\ifShort
We were able to achieve the current demand satisfaction in few minutes, 
while, today, managers need some hours.
Moreover, by setting a time limit of one hour,
we obtained solutions that
\fi
outperformed by far the ones of our customer, 
both in the missing demand and in the constraints violation 
\ifShort
(currently managers can violate some of them 
to have better schedules).
\fi

The paper is organized as follows.
In the rest of this section we mention the related work.
Section~\ref{sec:model} presents the mixed-integer linear programming
\ifShort 
formulation for the SDPTSP.
\fi
In Section~\ref{sec:results} we show some experimental 
\ifShort
results.
\fi
Finally, Section~\ref{sec:conclusion} concludes the paper and 
introduces the future research directions. 

\subsection{Literature review}\label{sec:lit}

The literature work related to personnel, staff, roster, 
or crew scheduling, is 
huge: there are thousands of papers 
studying several different variants of this problem. 
\ifShort
Due to the vast literature, we focus 
only on the few works with big importance for our study case, 
by mentioning the related problems
and suggesting some references.
\fi 

The problem of creating shifts is called Minimum Shift Design problem (MSD)
and various solution methods have been proposed, 
\ifShort
for more details see~\cite{digaspero}.
\fi
\ifShort
The assignment of employees to fixed shifts, called personnel scheduling,
has been introduced in the 1950s by Dantzig~\cite{dantzig}.
\fi 
Solution methods have been classified into several categories,
we refer to Alfares~\cite{alfares} for a comprehensive survey.
The problems concerned with assigning a set of tasks 
with fixed times to a heterogeneous workforce
having predetermined working times fall under the name of
Personnel Task Scheduling problems (PTSP), see~\cite{krishnamoorthy}.
Considering also the design of the shifts makes the situation harder 
and defines the Shift-Design Personnel Task Scheduling Problem (SDPTSP). 
As far as we know, only two papers dealt with this problem, 
they used as objective function an Equity criterion, 
for which they proposed two different approaches: 
one based on constraint programming~\cite{SDPTSP-E-CP}
and a two-step metaheuristic~\cite{SDPTSP-E-MH}.

The mentioned problems have lots of real-world applications,
Ernst et al.~\cite{ernst} presents some of the possible
\ifShort
applications areas.
\fi
For a very broad and detailed survey on personnel scheduling 
and rostering problems we refer to~\cite{reviewAll}, 
another interesting review combining managerial insights 
and technical knowledge is~\cite{reviewMan}.

\section{Definition of the MILP for the SDPTSP} \label{sec:model}

\ifShort
We consider the Shift-Design Personnel Scheduling Task Problem:
this means that given the two sets of the employees, 
with the individual skills and the time availability, 
and of the demands, with the associated task and times,
we want to design the shifts,
by fixing the times and the tasks,
and to assign them to the different employees.
The objective is to cover 
as much of the demand as possible
while respecting legal obligations 
(max daily hours, max daily span, ...)
and company constraints (min time on same activity, 
min work before break, ...).
\fi

\ifShort
We denote by $R$, $A$, $D$, and $T$ the sets
of employees, activities, days, and time slots, respectively.
\fi


\ifShort
\fi
The set of days is defined as
$D = \{ 1, \dots , \lvert D \rvert \}$, where $\lvert D \rvert$ 
is the planning horizon.
The set of time slots, denoted by $T$, 
contains numbers representing the starting minute. 
For example, if the interval length is 30 minutes, 
so the time step is $ts = 30$, 
we represent the time slot 8-8:30 
by the number $16 = (8 \cdot 60) / ts$, 
the time slot 8:30-9 by $17 = (8 \cdot 60 + 30) / ts$, ... ,
hence the set of time slots is $T = \{ 16, 17, \dots \}$.

\ifShort
To simplify the formulation,
\fi
we denote by $D_k$ the set of the first $k$ days, i.e.,
$D_k = \{ 1, \dots , k \}$, and its complementary,
so the set of all days but the first $k$,
by $D_{\overline{k}} = D \setminus D_k$. 
In the same way $D^k$ is the set of the last $k$ days,
and the set of all the days but the last $k$
is denoted by $D^{\overline{k}}$.
Similarly we define $T_k$ (resp. $T^k$) as the set of the first 
(resp. last) $k$ time slots and $T_{\overline{k}}$ (resp. $T^{\overline{k}}$) 
as the set of time slots but the first (resp. last) $k$.

Finally, we denote by $tsD$ the number of time slots in a day,
so $tsD = 1440/ts$.

\subsection{Variables}

The main variable that we use is the assignment of employee $r$ to 
activity $a$ at time slot $t$ of day $d$, 
\ifShort
so $x(r,a,t,d)=1$ if $r$ does $a$ at $t$ of $d$ and 0 otherwise. 
\fi
Another important information is the \textit{change of work status} 
of employee $r$ at time slot $t$ of day $d$ on activity $a$,
\ifShort
i.e., $y(r,a,t,d)=1$ (resp. $\unaryminus 1$) 
if $r$ starts (stops) $a$ at $t$ of $d$, while 
\fi
$y(r,a,t,d)=0$ means that $r$ on day $d$ either works both at $t$ and 
$t\unaryminus 1$
on $a$ or that he does not perform $a$ in any of these two intervals.
\ifShort
We also need a variable 
indicating whether employee $r$ works on day $d$, 
we have that $z(r,d)=1$ if $r$ works on at least one time slot of $d$ 
and $z(r,d)=0$ otherwise.
\fi
Moreover for each employee $r$ and for each day $d$ we define 
the begin (resp. end) of work,
which is represented by the first (resp. last) 
time slot with some activity assigned,
\ifShort
which is defined $b(r,d)$ (resp. $e(r,d)$).
\fi

The demands are defined by exploiting the activity $a$, the day $d$,
and the start and end time, namely $t_1$ and $t_2$ (such that $t_1 < t_2$).
It follows that $dem(a,d,t_1,t_2)=k$ means that 
on day $d$ between time slots $t_1$ and $t_2$
employees should work $k$ minutes on activity $a$. 
In order to keep the problem feasible, 
we add a nonnegative slack on each demand,
\ifShort
so $\alpha(a,d,t_1,t_2)$ are the minutes of demand $dem(a,d,t_1,t_2)$ 
which is not satisfied.
\fi

\subsection{Constraints}\label{sec:modelC}

\ifShort
In this section we present all the constraints of the model,
by exploiting their mathematical formulation.
They can be divided into three different groups:
problem defining, 
i.e., giving the relationships between variables,  
legal, i.e., fixed by the law, and 
company, i.e., defined by our customer 
to have better schedules. 
\fi
In order to guarantee feasibility 
the satisfaction of the demand is modeled as a soft constraint.
\ifShort
The compatibility between employees and activities depends 
on the individual skills of each worker
and are defined by the compatibility matrix, which is created 
in a pre-processing step by matching the skills of 
the employees and the requirement for each activity.
in this way we avoid to consider the skills in the model.
\fi

All the variables are derived from the assignment variable $x(r,a,t,d)$, 
the constraints defining these relationships are
\begin{enumerate}[label=(\enumWithLetter{D}{{\arabic*}})]
    \item definition of $y$
    \begin{align*}
    y(r,a,t,d) & = x(r,a,t,d) && 
    \quad \forall \ r \in R, a \in A, \ t \in T_1, \ d \in D \\
    y(r,a,t,d) & = x(r,a,t,d) - x(r,a,t \unaryminus 1,d) && 
    \quad \forall \ r \in R, a \in A, \ t \in T_{\overline{1}}, \ d \in D
    \end{align*}
    
    \item definition of $z$
    \begin{align*}
    \sum_{a \in A} \sum_{t \in T} x(r,a,t,d) & \le z(r,d) \cdot M && 
    \qquad \forall \ r \in R, \ d \in D \\
    \sum_{a \in A} \sum_{t \in T} x(r,a,t,d) & \ge z(r,d) && 
    \qquad \forall \ r \in R, \ d \in D
    \end{align*}
    where $M$ is a sufficiently big number, 
    e.g., $M = \lvert A \rvert \cdot \lvert T \rvert$,
    
    \item definition of $b$ and $e$
    \begin{align*}
    b(r,d) & \le t + tsD 
    \cdot \big(1 - \sum_{a \in A} x(r,a,t,d) \big) &&
    \qquad \forall \ r \in R, \ t \in T, \ d \in D \\
    e(r,d) & \ge (t+1) \cdot \sum_{a \in A} x(r,a,t,d) &&
    \qquad \forall \ r \in R, \ t \in T, \ d \in D
    \end{align*}
    
    \item each employee can perform simultaneously at most one activity
    $$ \sum_{a \in A} x(r,a,t,d) \le 1 \qquad \forall \ r \in R, 
    \ t \in T, \ d \in D$$

    \item correctness of start and end times, so
    $$ b(r,d) \le e(r,d) \qquad \forall \ r \in R, \ d \in D$$
\end{enumerate} 

\ifShort
The legal constraints depend on the work laws of the country,
for us they are
\fi
\begin{enumerate}[label=(\enumWithLetter{L}{{\arabic*}})]
    \item maximum daily hours
    \ifNoK 
    (\texttt{DH} minutes)
    $$ \sum_{a \in A}  \sum_{t \in T} x(r,a,t,d) \le \texttt{DH}/ts 
    \qquad \forall \ r \in R , \ d \in D$$
    \else ($k_{L1}$ time slots)
    $$ \sum_{a \in A}  \sum_{t \in T} x(r,a,t,d) \le k_{L1} 
    \qquad \forall \ r \in R , \ d \in D$$
    \fi
    
    \item maximum working hours in the planning horizon 
    \ifNoK 
    (\texttt{WH} minutes)
    $$\sum_{a \in A} \sum_{t \in T} \sum_{d \in D} x(r,a,t,d) \le \texttt{WH}/ts
    \qquad \forall \ r \in R$$
    \else ($k_{L2}$ time slots)
    $$\sum_{a \in A} \sum_{t \in T} \sum_{d \in D} x(r,a,t,d) \le k_{L2}
    \qquad \forall \ r \in R$$
    \fi
        
    \item 
    \ifNoK maximum \texttt{CD} consecutive working days 
    $$\sum_{d' = d}^{d+\texttt{CD}} z(r,d') \le \texttt{CD} 
    \qquad \forall \ r \in R, \ d \in D^{\overline{{\texttt{CD}}}}$$
    \else maximum $k_{L3}$ consecutive working days 
    \ifShort
    $$\sum_{d' \in [d, d+k_{L3}]} z(r,d') \le k_{L3} 
    \qquad \forall \ r \in R, \ d \in D^{\overline{{k_{L3}}}}$$
    \else
    $$\sum_{d' = d}^{d+k_{L3}} z(r,d') \le k_{L3} 
    \qquad \forall \ r \in R, \ d \in D^{\overline{{k_{L3}}}}$$
    \fi
    \fi  
    
    \item maximum consecutive working time, 
    \ifNoK 
    namely every \texttt{WT\textsubscript{1}} + \texttt{WT\textsubscript{2}} 
    minutes, each employee must be off for at least 
    \texttt{WT\textsubscript{2}} minutes, which is formulated
    $$ \sum_{t' = t}^{t+(\texttt{WT\textsubscript{1}}+\texttt{WT\textsubscript{2}})/ts \unaryminus 1} \sum_{a \in A} 
    x(r,a,t',d) \le \texttt{WT\textsubscript{1}}/ts 
    \qquad \forall \ r \in R, \ t \in T^{\overline{(\texttt{WT\textsubscript{1}}
    + \texttt{WT\textsubscript{2}})/ts \unaryminus 1}}$$
    \else 
    namely every $k_{L4} + j_{L4}$ time slots 
    each employee must be off for at least $j_{L4}$ slots,
    which is formulated
    \ifShort
    $$ \sum_{t' \in [t,t+k_{L4}+j_{L4})} \sum_{a \in A} 
    x(r,a,t',d) \le k_{L4} 
    \qquad \forall \ r \in R, \ t \in T^{\overline{k_{L4}+j_{L4} 
    \unaryminus 1}}$$
    \else
    $$ \sum_{t' = t}^{t+(k_{L4}+j_{L4} \unaryminus 1)} \sum_{a \in A} 
    x(r,a,t',d) \le k_{L4} 
    \qquad \forall \ r \in R, \ t \in T^{\overline{k_{L4}+j_{L4} \unaryminus 1}}$$
    \fi
    \fi
    
    
    \item maximum daily span 
    \ifNoK (\texttt{DS} minutes) 
    $$ e(r,d) - b(r,d) \le \texttt{DS} / ts \qquad \forall \ r \in R, \ d \in D$$
    \else ($k_{L5}$ time slots)
    $$ e(r,d) - b(r,d) \le k_{L5} \qquad \forall \ r \in R, \ d \in D$$
    \fi
    
    \item minimum rest between two working days 
    \ifNoK (\texttt{mr} minutes) 
    $$ tsD + b(r,d) - e(r, d \unaryminus 1) \ge \texttt{mr} / ts
    \qquad \forall \ r \in R, \ d \in D_{\overline{1}}$$
    \else ($k_{L6}$ time slots)
    $$ tsD + b(r,d) - e(r, d \unaryminus 1) \ge k_{L6} \qquad \forall \ r \in R, \ d \in D_{\overline{1}}$$
    \fi
\end{enumerate}

\ifShort
During several meetings with our customer we defined and implemented
their internal constraints,
i.e., the company constraints, which are
\fi
\begin{enumerate}[label=(\enumWithLetter{G}{{\arabic*}})]
    \item minimum working time after a break 
    \ifNoK (\texttt{wt} minutes) 
    $$ \sum_{t' = t}^{t+(\texttt{wt}/ts \unaryminus 1)} \sum_{a \in A} x(r,a,t',d) 
    \ge \sum_{a \in A} y(r,a,t,d) \cdot \texttt{wt} / ts
    \qquad \forall \ r \in R, \ t \in T^{\overline{\texttt{DM}/ts \unaryminus 1}}, \ d \in D$$
    \else ($k_{G1}$ time slots) 
    \ifShort
    $$ \sum_{t' \in [t,t+k_{G1})} \sum_{a \in A} x(r,a,t',d) 
    \ge \sum_{a \in A} y(r,a,t,d) \cdot k_{G1} 
    \qquad \forall \ r \in R, \ t \in T^{\overline{k_{G1} \unaryminus 1}}, \ d \in D$$
    \else
    $$ \sum_{t' = t}^{t+(k_{G1} \unaryminus 1)} x(r,a,t',d) \ge \sum_{a \in A} y(r,a,t,d) \cdot k_{G1} 
    \qquad \forall \ r \in R, \ t \in T^{\overline{k_{G1} \unaryminus 1}}, \ d \in D$$
    \fi
    \fi
    
    \item demand satisfaction, formulated as
    $$\sum_{t \in [t1,t2)} \sum_{r \in R} x(r,a,t,d) \cdot ts
    + \alpha(d,a,t_1,t_2)
    \ge dem(d,a,t_1,t_2) \quad \forall \ a \in A, \ d \in D, \ t_1,t_2 \in T$$

    \item minimum consecutive working time on the same activity 
    \ifNoK (\texttt{w\textsubscript{a}} minutes) 
    $$ \sum_{t' = t}^{t+(\texttt{w\textsubscript{a}}/ts \unaryminus 1)} x(r,a,t,d) 
    \ge y(r,a,t,d) \cdot \texttt{w\textsubscript{a}} / ts
    \qquad \forall \ r \in R, \ a \in A, \ t \in T^{\overline{k_a \unaryminus 1}}, \ d \in D$$    
    \else ($k_{G3,a}$ time slots) 
    \ifShort
    $$ \sum_{t' \in [t,t+k_{G3,a})} x(r,a,t,d) \ge y(r,a,t,d) \cdot k_{G3,a} 
    \qquad \forall \ r \in R, \ a \in A, \ t \in T^{\overline{k_{G3,a} \unaryminus 1}}, \ d \in D$$    
    \else
    $$ \sum_{t' = t}^{t+(k_{G3,a}-1)} x(r,a,t,d) \ge y(r,a,t,d) \cdot k_{G3,a} 
    \qquad \forall \ r \in R, \ a \in A, \ t \in T^{\overline{k_{G3,a} \unaryminus 1}}, \ d \in D$$    
    \fi
    \fi
    
    \item compatibility between employee and activity
    $$ x(r,a,t,d) \le cRA(r,a) 
    \qquad \forall \ r \in R, \ a \in A, \ t \in T, \ d \in D$$
    where $cRA$ is the compatibility matrix matching employees and activities
    such that $cRA(r,a)=1$ if $r$ can do $a$ and 0 otherwise
    
    \item compatibility between employee and time slot
    $$ x(r,a,t,d) \le cRTD(r,t,d)
    \qquad \forall \ r \in R, \ a \in A, \ t \in T, \ d \in D$$
    where $cRTD$ is the availability matrix such that $cRTD(r,t,d) = 1$
    if $r$ is available to work at $t$ of $d$ 
    and 0 otherwise
    
    \item the checkout management comprises two different activities,
    i.e., normal working \textit{opCAS} and closure \textit{clCAS}, 
    and it has a specific daily rule:
    ``each employee who works at the checkout 
    has to do only one closure rigth after its last opening'':
    so, for each day, if an employee works on the checkout 
    (do at least one slot \textit{opCAS}),
    then he must do one slot \textit{clCAS} 
    right after its last \textit{opCAS} slot; 
    this is modelled by the following constraints
    \begin{itemize}
        \item an employee can do \textit{clCAS} only if, 
        the same day, he works also on \textit{opCAS} 
        $$ \sum_{t \in T} x(r,opCAS,t,d) \le \sum_{t \in T} x(r,clCAS,t,d) \cdot M
        \qquad \forall \ r \in R, \ d \in D$$
        for $M$ sufficiently big, e.g., $M = \lvert T \rvert $

        \item each employee does at most one \textit{clCAS} per day
        $$ \sum_{t \in T} x(r,clCAS,t,d) \le 1 \qquad \forall \ r \in R, \ d \in D$$
    
        \item after doing \textit{clCAS}, 
        the employee cannot do \textit{opCAS} until the next day
        $$ \sum_{t' \ge t} x(r,opCAS,t',d) \le M(1 - x(r,clCAS,t,d))
        \qquad \forall \ r \in R, \ t \in T, \ d \in D$$
        for $M$ sufficiently big
        
        \item right before \textit{clCAS} the employee must do 
        a \textit{opCAS} time slot
        $$ x(r,clCAS,t,d) \le x(r,opCAS,t \unaryminus 1,d) 
        \qquad \forall \ r \in R, \ t \in T_{\overline{1}}, \ d \in D$$
    \end{itemize}    

\end{enumerate}

\ifShort
It is important to note that some constraints, e.g., $(L6)$, 
have to consider the previous schedules;
in the results of the next section we also add the historical information
even though we did not report the exact formulation 
of these constraints throughout this section.
\fi

\subsection{Objective function}\label{sec:modelOF}

\ifShort
In the optimization process several objective functions might be used,
in general the aim is to minimize the costs 
(declined either as total working days or as uncovered demand).
\fi

\ifShort
The company we are working with asked us, as objective function, 
to cover the most demand possible
in order to avoid the employees doing overwork.
They also specified that the activities have different priorities, 
so we gave a penalty to the slack for each activity $a$, denoted $p(a)$. 
Moreover, their assignment of employees to tasks 
depends on the individual work experience,
so an employee should be assigned to most suited activity even
if he can perform several of them. 
In order to catch this, 
we added a multiplicative factor which depends on the matching
between employee $r$ and activity $a$, denoted $c(r,a)$.
\ifShort
\else
Since we are minimizing the value of the objective function a smaller cost 
corresponds to a better matching, 
in other words $c(r,a_1) < c(r,a_2)$ means that employee $r$
would rather work on $a_1$ than on $a_2$.
\fi
Thus, the objective function used in the model is
$$ \min \sum_{a \in A} \sum_{d \in D} \sum_{t_1, t_2 \in T} 
p(a) \cdot \alpha(a,d,t_1,t_2) + \sum_{r \in R} \sum_{a \in A} 
\sum_{t \in T} \sum_{d \in D} x(r,a,t,d) \cdot c(r,a)$$

\fi

\section{Instances and Computational Aspects} \label{sec:results}

The main scope of this paper is to introduce 
the mixed-integer linear programming formulation
for the Shift-Design Personnel Task Scheduling Problem,
with the specific constraints for our case study.
\ifShort
In order to validate this formulation 
and to provide relevance to this problem, 
\fi
we solve three benchmark instances 
corresponding to three different planning weeks of our customer
with a granularity of 15 minutes
(so $\lvert T \rvert = 48$, and $\lvert D \rvert = 7$),
and we compare the solution with the actual schedules.
The instances have small-medium size, because
$\lvert R \rvert = 66$, and $\lvert A \rvert = 72$.

The following results are obtained by using 
the commercial solver Xpress~\cite{fico}, 
with version 8.11 on a standard laptop
Intel Core i7-8550 at 1.8 GHz with 16 GB RAM running Windows 10.

\subsection{Resolution Steps}\label{sec:results1}

During the analysis phase and the study of the first results
we were able to find some interesting considerations 
that led us to reduce the size of the problem.
\ifShort
From a computational point of view, the most useful
was the \textit{identity} of some activities, i.e.,
there are some tasks 
which can be performed by the same set of employees, 
with nonoverlapping demand, and for which all the parameters are equal.
Hence, we merged the identical tasks into macro-activities 
to reduce the size of the problem.
\fi
By doing this as a pre-processing step 
we were able to simplify the problem
and to obtain better results:
we decreased the number of activities, on average, by 29\%, 
and we improved the best solution by 37\%; 
the detailed results can be seen in Table~\ref{table:optitRes}.
\ifShort
After the optimization, the assignment to these macro-activities
is redistributed on the original tasks by a simple backward assignment procedure.
\fi

As it is easy to imagine, 
the main difficulty of a commercial solver
on these kind of problems is to give big improvements on the best solution: 
since there are many feasible solutions with similar values, 
the incumbent continuously decreases but very slowly.
\ifShort
We solved this problem by implementing a base heuristic algorithm, 
it produced a feasible assignment in few seconds
which was then passed to the solver as a warm-start.
The heuristic we developed uses a greedy approach: 
it creates the shifts for all the employees satisfying the compatibilities
of tasks and hours,
then it checks for all the conflicts on the constraints 
and it fixes them by changing the tasks, 
in this phase the priority is given to activities with higher demand.
\fi

\subsection{Experimental Results}\label{sec:results2}

In this section we present the results 
obtained on a set of three instances provided 
by the company we are working for, 
which correspond to three week of plannings 
of a medium-big size retail store.

\ifShort
Since the results we obtained were already rewarding, 
we agreed with our customer to set a time limit of one hour,
which is far less the actual time needed (around 6 hours) 
by the store managers.
\fi

\begin{table}[ht]
\centering
\begin{tabular}{|c|c|c|c|c|c|}
\hline
Instance & Merge & Heur & Best sol & Best bound & Gap \\ 
\hline\hline
W20 & N & N & 701 384 & 23 073 & 97\% \\ 
\hline
W20 & Y & N & 369 069 & 23 041 & 94\% 
\\\hline
W20 & N & Y & 87 691 & 23 073 & 74\% \\ 
\hline
W20 & Y & Y & 95 517 & 23 041 & 76\% \\ 
\hline \hline
W21 & N & N & 556 049 & 4 780 & 99\% \\
\hline
W21 & Y & N & 459 375 & 4 778 & 99\% \\
\hline
W21 & N & Y & 103 632 & 4780 & 95\% \\
\hline
W21 & Y & Y & 77 633 & 4 778 & 94\% \\
\hline \hline
W22 & N & N & 693 821 & 31 009 & 95\% \\ 
\hline
W22 & Y & N & 297 904 & 30 976 & 89\% \\
\hline
W22 & N & Y & 169 001 & 31 009 & 82\% \\
\hline
W22 & Y & Y & 159 152 & 30 976 & 80\% \\
\hline
\end{tabular}
\caption{Computational results of Optit 
for the instances given by our customer 
with a time limit of one hour}
\label{table:optitRes}
\end{table}

In Table~\ref{table:optitRes} we validated 
the computational approaches presented in Section~\ref{sec:results1}
by exploiting the results for the different combinations 
of the pre-processing techniques explained.
The second column indicates whether \textit{identical} activities 
have been merged,
while the third if the greedy heuristic has been used as a warm start.
From these results it follows that using the heuristic 
decreases the best solution on average by 78\%, 
while merging the activities gives an improvement of 37\%.
The total processing time is around 1 minute, 
so it is negligible compared to the time limit.

We claim that these results are far from being optimal,
and it can be seen from the gaps in Table~\ref{table:optitRes},
this is mainly because our approach is not problem-dependent:
we used a commercial solver on a general formulation
without tuning the internal parameters 
and we did not implement specific heuristics. 
Anyway, these first results are very promising for us and 
rewarding for our customer.

Since these instances are provided by the company we are working with, 
there is no benchmark result to check the quality of the solution.
The only comparison we can make 
\ifShort
is with the current plannings of our customer.
\fi

\begin{table}[ht]
\centering
\begin{tabular}{|c|c|c|c|c|c|}
\hline
Instance & Case & $\sum \alpha$ [h] & Violations & Dept dem \% \\ 
\hline\hline
W20 & C & 374 & 579 & 86\% \\ 
\hline
W20 & O & 59 & 0 & 96\% \\
\hline\hline
W21 & C & 366 & 544 & 86\% \\ 
\hline
W21 & O & 62 & 0 & 95\% \\
\hline\hline
W22 & C & 356 & 744 & 86\% \\ 
\hline
W22 & O & 110 & 0 & 95\% \\
\hline
\end{tabular}
\caption{Comparison of Optit results (case $O$) 
with actual planning (case $C$)}
\label{table:optitDimar}
\end{table}

In Table~\ref{table:optitDimar} we compare our results 
(case $O$) with the real schedules of our customer (case $C$).
In the third column we provide the total missing hours 
to cover the demand, i.e., the sum of the slacks 
without any multiplicative parameter,
while the fourth indicates the number of violated constraints;
since we use an exact model, 
we do not violate any constraint,
but in the real life managers do it quite often.
In the fifth column we give a percentage of saturation of the whole 
demand by considering the departments and not the single tasks, 
this is done because sometimes the managers 
assign the employees to only one task, 
but then they can freely shift activity
within the same department if needed.

It is easy to see that the missing demand in our solution
is much smaller than the actual planning.
If we consider the single activities we have a huge decrease 
of around 70\%.
\ifShort
With respect to the associated departments, 
we can see that on average our demand satisfaction is 10\% more, 
which still gives a big improvement.
\fi
\ifShort
Another important indicator of the quality of our solutions 
are the constraints violated in the current weekly schedules (around 600),
which for us are 0 since we use an exact model.
\fi

It follows from the results that, 
in comparison with the actual planning, 
we are outperforming the current schedules
both with respect to the covered demand 
and to the violation of the constraints.

\section{Conclusion}\label{sec:conclusion}

\ifShort 
The Shift-Design Personnel Task Scheduling Problem has been studied 
only by few authors in the literature
and, to the best of our knowledge, 
no exact methods were proposed to solve this problem.
The main work of this paper has been to give the first MILP
formulation for this problem, 
defining the generic constraints given by the law.
Our study case is a big size Italian retail company,  
so we also implemented some constraints for this specific case,
to provide them with better results.
In order to validate our model 
we studied three different weeks of planning on a store, 
and we provided the results, by using the commercial solver Xpress.
Moreover, we explained how to decrease the magnitude of the problem, 
in order to tackle more difficult, i.e., bigger, instances.
\fi
Finally, we compared our results with the actual planning.
We were able to decrease the uncovered demand by 70\%, 
which corresponds to a huge saving.
At the same time our solutions respected all the constraints, 
while, as of today, the schedules of our customer violate 
around 600 constraints every week.
Considering the department and not the single tasks, the improvement is lower,
but still very significant since we cover 10\% more demand.
It is also important to say that as of today managers need six hours every week
for the planning, while we set a time limit of 60 minutes.
Hence we obtained better results in a less amount of time.

\ifShort
As it can be imagined the company we are working with 
is very satisfied with these first results, 
but there is still ongoing work in different directions.
As a first step, we are trying to catch some other needs, 
in order to add more specific constraints 
and to provide even better results.
Another direction of future work is the development
of specific techniques to speed up the computational times 
and, hence, to solve bigger instances. 
We saw that, 
by passing a feasible schedule to the solver as a warm-start, 
the quality of the solution improved drastically 
in a short amount of time,
even with our simple greedy heuristic: 
we were able to reach the same results 
in few minutes instead of an hour.
Thus, implementing a smarter heuristic 
is one of the tasks with higher priority
to further improve the current solutions.
\fi

%

%
%

\end{document}